\documentclass{amsart}

\usepackage{amsfonts,amssymb,amsmath,amscd,amsthm}
\usepackage[mathscr]{eucal}
\usepackage{mathrsfs}

\usepackage{mathbbol}

\usepackage{oldgerm,units}
\usepackage{wrapfig,epsfig}

\usepackage{ifthen}
\newtheorem{theorem}{Theorem}[section]
\newtheorem{proposition}[theorem]{Proposition}
\newtheorem{definition}[theorem]{Definition}

\newtheorem{lemma}[theorem]{Lemma}

\newtheorem{corollary}[theorem]{Corollary}

\newtheorem{example}[theorem]{Example}
\newtheorem{remark}[theorem]{Remark}

\newcommand{\Real}{\mathbb R}

\newcommand{\one}{\mathbb{1}}
\newcommand{\zero}{\mathbb{0}}

\newcommand{\Trop}{\mathbb T}
\newcommand{\RealInf}{{\Real \cup \{ \tUniS \}}}

\newcommand{\tUnit}{\mathbb U}
\newcommand{\etUnit}{\bar{\tUnit}}
\newcommand{\eReal}{\bar{\Real}}

\newcommand{\tA}{A}
\newcommand{\tB}{B}

\newcommand{\tI}{I}

\newcommand{\tV}{v}

\newcommand{\tX}{x}

\newcommand{\To}{\longrightarrow }

\newcommand{\Dir}{\thSpc \Longrightarrow \thSpc}

\newcommand{\dir}{\Rightarrow}
\newcommand{\tUniS}{-\infty}

\newcommand{\uuu}[1]{#1^\nu}

\newcommand{\epiToMaxPlus}{\pi}

\newcommand{\invA}[1]{#1^{\bigtriangledown}}% {#1^{\bigtriangledown}}
\newcommand{\al}{\alpha}
\newcommand{\bt}{\beta}
\newcommand{\gm}{\gamma}
\newcommand{\sig}{\sigma}

\newcommand{\lm}{\lambda}
\newcommand{\Lm}{\Lambda}

\newcommand{\itA}{\invA{\tA}}

\newcommand{\Genr}[1]{\widetilde{#1}}
\newcommand{\gtI}{\Genr{\tI}}

\newcommand{\rnk}{\operatorname{rk}}

\newcommand{\TrS}{\oplus}
\newcommand{\TrP}{\odot}

\newcommand{\TropSR}{(\Trop,\TrS,\TrP)}

\newcommand{\eMaxPlusAlg}{(\eReal,\max,+)}

\newcommand{\OP}{\left(}
\newcommand{\CP}{\right)}

\def \aaa{a}
\def \aab{b}

\def \oDeg{d_{\operatorname{out}}}
\def \iDeg{d_{\operatorname{in}}}

    \ifx\proof\undefined
    \newenvironment{proof}{
    \smallskip
    \noindent\emph{Proof.}}{\hfill\(\Box\)
    \bigskip
    } \fi

\newcommand{\vvMat}[9]{\small{\OP \begin{array}{ccc}
  #1 & #2 & #3\\
  #4 & #5 & #6\\
  #7 & #8 & #9\\
\end{array}\CP}}

\newcommand{\bfem}[1]{\textbf{\emph{#1}}}

\newcommand{\ifdef}[3]{\ifthenelse{\equal{#1}{true}}{#2}{#3}}

\newcommand{\thSpc}{\; \; \;}

\newcommand {\secSpc} {\vskip 0cm}

\newcommand {\parSpc} {\vskip 0cm}

\def\trn{{\operatorname{t}}}
\newcommand{\Adj}[1]{\operatorname{adj}({#1})}
\def\Size{\operatorname{size}}

\newcommand{\etype}[1]{\renewcommand{\labelenumi}{(#1{enumi})}}
\def\eroman{\etype{\roman}}

\def\pSkip{\vskip 1.5mm \noindent}
\def\piSkip{\vskip 1.5mm }

\def\({\left(}
\def\){\right)}

\newcommand{\si}[2]{{#1, #2}}

\def\bal{\bar \al}
\def\bbt{\bar \bt}

\def\one{\mathbb{1}}
\def\zero{\mathbb{0}}

\def\msig{\hat \sig}

\newcommand{\matS}[2]{M_{#1 \times #2}}
\newcommand{\Det}[1]{ \left|{#1}\right|}

\def\cyc{C}

\def\pth{P}
\def\grph{G}
\def\out{{\operatorname{out}}}
\def\inn{{\operatorname{in}}}
\def\odeg{d_\out}
\def\ideg{d_\inn}

\def\nmulti{$n$-multicycle}

\def\minf{-\infty}

\def\bfa{{\bf a}}

\begin{document}

%******************************* title ***********************************

\title[The Tropical Rank of a Tropical Matrix]
{The Tropical Rank of a Tropical Matrix}

%******************************* authors *********************************

\author{Zur Izhakian}\thanks{The author has been supported by The German-Israeli Foundation for Research and
Development by Hermann Minkowski Minerva Center for Geometry at
Tel Aviv University.}
\address{Department of Mathematics, Bar-Ilan University, Ramat-Gan 52900,
Israel}  \address{ \vskip -6mm CNRS et Universit´e Denis Diderot
(Paris 7), 175, rue du Chevaleret 75013 Paris, France}
\email{zzur@math.biu.ac.il}

%******************************* AMS classification ***********************

%******************************* date *************************************
\date{July 2008}

%******************************* keywords *********************************


%******************************* abstract *********************************

%******************************* abstract *********************************
\begin{abstract}
In this paper we further develop  the theory of matrices over the
extended tropical semiring. Introducing a notion of tropical
linear dependence allows for a natural definition of matrix rank
in a sense that coincides with the notions of tropical regularity
and invertibility.
\end{abstract}

\maketitle

%******************************* remarks *********************************

%\tableofcontents

%******************************* section ******************************
%******************************* BODY *********************************

%******************************* SECT 1 *********************************
\secSpc
\section*{Introduction}\label{sec:Introduction}

One of the most important notions in linear algebra is the notion
of rank, especially with a suitable relation to linear dependence.
In the familiar tropical linear algebra the  notion of dependence
is absent, mostly since the ground max-plus semiring is
idempotent. The special structure of the \bfem{extended tropical
semiring}, as introduced in \cite{zur05TropicalAlgebra}, allows a
natural definition for this absent notion, providing the tropical
analogous to rank of matrices as in the classical theory, i.e. the
maximal number of independent rows, and leading to the two
important results:
\begin{itemize}
    \item An $n \times n$  matrix $\tA$ has rank $n$ iff $\tA$ is tropically
    nonsingular iff $A$ is pseudo invertible, \pSkip
    \item An $m \times n$ matrix $\tA$ has rank $k$ iff its maximal
nonsingular minor is of size $k \times k$.

\end{itemize}
Although our framework is typically combinatorial, as these
results show, the tropical analogous to classical results are
carried naturally over the extended tropical semiring.

\parSpc
The main goal of this paper is a further development of
 the basics of tropical matrix algebra over the
{extended tropical semiring}, $\TropSR$, as has been presented in
\cite{zur05TropicalAlgebra}; we also use some of the terminology
used in \cite{IzhakianRowen2007SuperTropical}.
 This extension is obtained by taking two copies of the reals,
$$\eReal = \RealInf \quad \text{ and } \quad \etUnit = \uuu{\Real} \cup \{ \tUniS
\},$$ each is enlarged by $\{\tUniS\}$, and gluing them along
$\tUniS$ to define the set $\Trop = \eReal \cup \etUnit$.  We
define the correspondence $\nu : \Real \to \tUnit$ to be the
identity map, and denote the image of $a \in \Real$ by $\uuu{a}$.
Accordingly, elements of $\tUnit$, which is  called the
\bfem{ghost} part of $\Trop$, are denoted as $\uuu{\aaa}$; $\Real$
is called the \bfem{real} part of $\Trop$. The map $\nu$ is
sometimes extended to whole $\Trop$,
\begin{equation}\label{eq:ghostMap}
    \nu: \Trop \To \etUnit,
\end{equation}
by declaring $\nu: \aaa^\nu \mapsto \aaa^\nu$ and $\nu : \tUniS
\mapsto \tUniS$. (We use the generic notation
 $\aaa, \aab\in \Real$ and $x ,y  \in \Trop$.)

The set $\Trop$ is then provided with the following total order
extending the usual order on $\Real$: \piSkip
\begin{enumerate}\eroman
    \item $\tUniS \prec x ,$ $\forall x  \in \Trop$; \piSkip
    \item for any real numbers $\aaa < \aab$, we have $\aaa \prec
\aab,$ $\aaa \prec \uuu{\aab}$, $\uuu{\aaa} \prec \aab$, and
$\uuu{\aaa} \prec \uuu{\aab}$; \piSkip
    \item  $\aaa \prec \uuu{a}$ for all $\aaa \in \Real$.
\end{enumerate}
  Then, $\Trop$ is endowed with the two operations $\TrS$
and $\TrP$ , defined as follows: \pSkip
\ \ $
\begin{array}{rll}
  x  \TrS y  & = &\left\{ \begin{array}{ll}
  max_{(\prec)}\{x , y \}, & x  \neq y , \\[2mm]
  \uuu{x }, & x  = y  \neq \tUniS,
\end{array}\right. \\[2mm] \tUniS \TrS \tUniS & = & \tUniS,
\end{array}
$
$
\begin{array}{rll} \aaa \TrP \aab & = & \aaa + \aab,\\[2mm]
a^\nu \TrP b & = & a \TrP b ^\nu = a^\nu \TrP b^\nu = (a+b)^\nu, \\[2mm]
 (\tUniS) \TrP x   & = & x  \TrP (\tUniS) = \tUniS.
\end{array}$\pSkip
We usually write $x   y  $ for the product $x  \TrP y $,  for
short. Similarly, the division is written  $\frac{x }{y }$, for
$y  \neq \minf$.

 Triple $\TropSR$ is called the {extended tropical
semiring}; this semiring is nonidempotent commutative semiring,
since $\aaa \TrS \aaa = \uuu{\aaa}$, with the unit element
$\one_{\Trop} := 0$ and the zero element $\zero_{\Trop} :=
\tUniS$.
 This, and the fact that $(\Real, \TrP)$ is a
group and $(\etUnit, \TrS,\TrP)$ is an ideal, provides $\Trop$
with a more richer structure to which much of the theory of
commutative algebra can be transferred.

The connection with the standard tropical (max-plus) semiring is
established by the natural semiring epimorphism,
\begin{equation}\label{eq:epiTopicalSemiRings}
    \epiToMaxPlus: \TropSR \; \To \;  (\Real \cup \{
\tUniS\}, \max, + \;), \\
 \end{equation}
where $\epiToMaxPlus:  \uuu{\aaa} \mapsto  \aaa$, $\epiToMaxPlus:
\aaa \mapsto \aaa$ for all $\aaa \in \Real$, and $\epiToMaxPlus:
\tUniS \mapsto \tUniS$. (We write $\epiToMaxPlus(x )$ for the
image of $x  \in \Trop$ in $\eReal = \Real \cup \{ \tUniS\}$.)
This epimorphism induces epimorphisms of polynomial semirings,
Laurent polynomial semirings, and tropical matrices.

%************************* Acknowledgement ****************************
\parSpc
\textbf{\emph{ Acknowledgement}}:  The author would also like to
thank \emph{Prof. Eugenii Shustin} for our fertile discussions.

%\end{AMS}

%\parbox[position][height][inner-pos]{width}{text}
%\framebox[width][position]{text}

%~~~~~~~~~~~~~~~~~~~~~~~~~~~ SUB-SECT ~~~~~~~~~~~~~~~~~~~~~~~~~~~~~~~~~~~
\secSpc
\section{Tropical vector spaces}\label{sec:TropicalDependence}

As in the classical ring theory, the tropical space $\Trop^{(n)}$,
consisting of all  $n$-tuples $(x_1, \dots,x_n)$ with entries $x_i
\in \Trop$, is treated as a semiring module with addition, and
multiplication by $\al \in \Trop$, defined with respect to
$(\Trop, \TrS. \TrP)$. An $n$-tuple  $(x_1, \dots,x_n)\in
\Trop^{(n)}$ is called vector, and a vector having only ghost, or
$\tUniS$, entries (i.e. $(x_1, \dots,x_n)\in \etUnit^{(n)}$) is
termed \bfem{ghost vector}.

\begin{definition} A set $W = \{e_1, \dots, e_n\}\subset \Trop^{(n)}$ is a \textbf{classical
base} of $ \Trop^{(n)}$, if every element of $\Trop^{(n)}$ can be
written uniquely in the form $\bigoplus_{i=1}^n \al_i e_i$, where
$\al_i \in \Trop$.
\end{definition}

The \bfem{standard base} of $\Trop^{(n)}$ is defined as
$$e_1 = (0,\tUniS, \dots, \tUniS),
\quad e_2 = (\tUniS,0,\tUniS, \dots, \tUniS), \quad \dots, \quad
e_n = (\tUniS,\tUniS, \dots, 0).$$

\begin{definition}\label{def:tropicalDependent} A collection of
vectors  $v_1,\dots,v_m$ is said to be \textbf{tropically
dependent}
 if there exist $\al_1,\dots,\al_n \in \eReal$, but not all of them $\tUniS$,
 for which
$$
   \al_1 v_1 \TrS \cdots \TrS \al_m  v_m \in
   \etUnit^{(n)},
$$
otherwise the vectors are said to be \textbf{tropically
independent}. We call these $\al_i$'s the \textbf{dependence
coefficients} of $v_1,\dots,v_m$.
\end{definition}
Any set of vectors containing a ghost vector is tropically
dependent; in particular, a singleton consisting of  a ghost
vector is tropically dependent.

\begin{example} Let $v_1,v_2,v_3$, and $v_4$  be the following tropical vectors:
$$ v_1 = (0,1), \qquad  v_2 = (1,2), \qquad \ v_3 = (2,0) , \qquad \ v_4 = (2^\nu,0).$$
Then, $v_1$ and $v_2$ are tropically dependent, since  $1 v_1 \TrS
v_2 \in \etUnit^{(2)}$. The vectors $v_1$ and $v_3$ are tropically
independent, but $v_1$ and $v_4$ are tropically dependent, i.e. $1
v_1 \TrS v_4 \in \etUnit^{(2)}$.
\end{example}

Different from the classical theory, where the ground structure is
a field, in which the notions of linear dependence and span
coincide, these notions do not coincide in the tropical framework.
Namely, even if a collection of vectors is linearly dependent it
might happen that no one can be expressed in terms of other
vectors; for example, take
$$ v_1 = (1,1,\tUniS), \quad   v_2 = (1,\tUniS,1), \quad \text{and
} \
 v_3 = (\tUniS,1,1), $$ these vectors are linearly dependent, i.e.
 $v_1 \TrS v_2 \TrS v_3 \in \etUnit^{(3)}$, but non of these vectors can be written
 in terms of the others.

%******************************* SECT  *********************************
\secSpc
\section{Regularity of tropical matrices}
%******************************* SECT  *********************************
%~~~~~~~~~~~~~~~~~~~~~~~~~~~ SUB-SECT ~~~~~~~~~~~~~~~~~~~~~~~~~~~~~~~~~~~
%\subSecSpc
\subsection{Tropical matrices}\label{sec:FundamentalofTropicalMatrices}

It is standard that since $\Trop$ is a semiring then we have the
semiring $\matS{n}{n}(\Trop)$ of $n \times n$ matrices with
entries in $\Trop$, where addition and multiplication are induced
from $\Trop$ as in the familiar matrix construction. The
\bfem{unit} element $I$ of
 $\matS{n}{n}(\Trop)$, is the matrix
with $0$ on the main diagonal and whose off-diagonal entries are
$\tUniS$; the \bfem{zero} matrix is  $Z = (\tUniS) \tI$;
therefore, $\matS{n}{n}(\Trop)$ is also a multiplicative monied.

We write $\tA =(a_{i,j})$ for a tropical matrix in
$\matS{n}{n}(\Trop)$ and denote the entries of $\tA$ as $a_{i,j}$.
We say that $A$ is \textbf{real matrix} if each $a_{i,j}$ is in
$\eReal$, $A$ is called \textbf{ghost matrix} when  each $a_{i,j}$
belong ot $\etUnit$. Since $\Trop$ is a commutative semiring, $ x
\tA = \tA x$ for any $x \in \Trop$. (We denote the set of $m
\times n $ matrices by $\matS{m}{n}(\Trop)$.)

As in the familiar  way, we define the \bfem{transpose} of $\tA =
(a_{i,j})$ to be $\tA^\trn = (a_{j,i})$. The \bfem{minor}
$\tA_{i,j}$ is obtained by deleting the $i$ row and $j$ column of
$\tA$. We define the \bfem{tropical determinant} to be
  \begin{equation}\label{def:tropicalDet}
 |\tA| = \bigoplus_{\sig \in S_n}  \OP  \aaa_{1, \sig(1)}
 \cdots \aaa_{n,\sig(n)} \CP,
\end{equation}
where $S_n$ is the set of all the permutations on $\{1,\dots,
  n\}$. Equivalently, the tropical determinant
$|\tA|$ can be written also in terms of minors as
\begin{equation}\label{def:tropicalDetByMinors}
|\tA| = \bigoplus_{j}   \aaa_{i ,j} |\tA_{i, j}|,
\end{equation}
for some fixed index $i$. Indeed, in the classical terminology,
since parity of indices' sums are not involved, the tropical
determinant is a permanent, what makes the tropical determinant a
pure combinatorial function.

 We use the notation $\msig$ for a permutation, not necessarily unique, whose
 $\nu$-evaluation in $A$
  equals $|A|^\nu$, and write
$$ \gm =  \aaa_{1,\msig(1)}
 \cdots \aaa_{n,\msig(n)};   $$
therefore $\pi\(\gm\) = \pi\(|A|\)$, or equivalently $\gm^\nu =
(|A|)^\nu$.  (We use both forms for convenience.)   We say that
two permutations $\sig_1$ and $\sig_2$ in $S_n$ are
\textbf{disjoint} if $\sig_1(i) \neq \sig_2(i)$ for each $i =
1,\dots,n$.

\begin{remark}\label{obs:detProp}
The tropical determinant has the following properties:
\begin{enumerate}
   \item Transposition and reordering of rows or columns leave the determinant
   unchanged; \pSkip
  \item The determinant is linear with respect to scalar multiplication of any given row or
  column by a real.
\end{enumerate}
\end{remark}

The \bfem{adjoint} matrix $\Adj{\tA}$ of a matrix $\tA
=(a_\si{i}{j})$ is defined as the matrix
$(a'_\si{i}{j})^{\operatorname{t}}$ where $a'_\si{i}{j}=
|\tA_\si{j}{i}|$.

  \begin{definition}\label{def:singularMat} A matrix $\tA \in \matS{n}{n}(\Trop)$
  is said to be \textbf{tropically singular}, or singular, for short,
  whenever
  $|\tA| \in \etUnit$, otherwise $\tA$ is called \textbf{tropically nonsingular},
  or nonsingular, for short.
  \end{definition}
  %
%\noindent
 In particular, when two or more different
permutations, $\msig_1, \msig_2, \dots  \in S_n$, achieve the
$\nu$-value of $|\tA|$ simultaneously, or the permutation $\msig$
that reaches the $\nu$-value of $|\tA|$ involves an entry in
$\etUnit$, then $\tA$ is singular.

\begin{remark}\label{rmk:regularity} In this combinatorial view, for real matrices, our definition of
singularity coincides with the known definition for matrices over
$\eMaxPlusAlg$, cf.~\cite{Develin2003}.
\end{remark}

To establish  a notion of pseudo invariability for
$\matS{n}{n}(\Trop)$, viewed as monoid, we define a \bfem{pseudo
unit matrix} to be a \textbf{regular} matrix with $0$ on the main
diagonal and whose off-diagonal entries are in $\etUnit$; in
particular, the unit matrix is also a pseudo unit. We use these
pseudo unit matrices to define the distinguished subset $U_{n
\times n} (\Trop) \subset \matS{n}{n}(\Trop)$ as
\begin{equation}\label{eq:setOfPartialUnits}
  U_{n \times n}(\Trop) = \left\{\gtI \; : \; \gtI \; \;  \text{ is a  pseudo unit
  matrix}
  \right\} \ .
\end{equation}
\begin{definition}\label{def:inverseMatrix}
A matrix $\tA \in \matS{n}{n}(\Trop)$ is said to be \textbf{pseudo
invertible}  if there exits a matrix $\tB \in \matS{n}{n}(\Trop)$
such that $\tA \tB \in U_{n \times n}(\Trop)$ and $\tB \tA \in
U_{n \times n}(\Trop)$. If $\tA$ is pseudo invertible, then we
call $\tB$  a \textbf{pseudo inverse matrix} of $\tA$ and denote
it as $\itA$.
\end{definition}
Having this setting, we  state one of our main theorems analogues
to the classical relation,
\cite[Theorem~3.3]{zur05TropicalAlgebra}:
\begin{theorem}\label{thm:inverseMatrix}
A matrix $\tA \in M_{n \times n}(\Trop)$ is pseudo invertible iff
is tropically regular. In case $\tA$ is regular, $\itA$ can be
defined as
$\itA = \frac{\Adj{\tA}}{|\tA|}$, and is called  the
\textbf{canonical pseudo inverse} of $\tA$.
\end{theorem}

%******************************* SECT  *********************************

%\subsection{The digraph of a tropical matrix}

%******************************* SECT  *********************************
\secSpc
\subsection{Lemmas on tropical regularity}

Our main computational tool in tropical matrix theory is the
\textbf{weighted digraph} $\grph = (V,E)$ of an $n\times n$ matrix
$A = (a_{i,j})$, which is defined to have vertex set $V =\{ 1,
\dots, n\}$, and an edge $(i,j)$ from $i$ to $j$ (given
\textbf{weight} $a_{i,j}$) whenever $a_{i,j} \ne  \minf$. We
denote this graph by $G_A$. In this view, reordering of rows or
columns of $A$ is equivalent to relabeling of vertices on $G_A$.

We use \cite{Gibbons85} as a general reference for graphs. We
always assume that $V = \{ v_1, \dots, v_n \}$, for convenience of
notation. The \textbf{out-degree}, $\odeg(v)$, of a vertex is the
number of edges emanating from $v$, and the \textbf{in-degree},
$\ideg(v)$, is the number edges terminating at $v$.  A
\textbf{sink} is a vertex with $\odeg(v) = 0$, while a
\textbf{source} is a vertex with $\ideg(v) = 0$.

 The \textbf{weight} $w(\pth)$ of a path $\pth$
is defined to be the sum of the weights of the edges comprising
$\pth$, counting multiplicity.
 A \textbf{simple cycle} is a simple path for
which $\odeg(v) = \ideg(v) =1$ for every vertex $v$ of the path;
thus, the initial and terminal vertices are the same.  A simple
cycle of length 1 is then a loop. We define a
$k$-\textbf{multicycle} $\cyc$ in a digraph to be the union of
disjoint simple cycles, the sum of whose lengths is $k$.

Writing a permutation $\sig$ as a product $\mu _1\cdots \mu_t$ of
disjoint cyclic permutations, we see that each permutation $\sig$
corresponds to an \nmulti\ in $G_A$, and their highest weight
matches~$|A|$. In particular, when $|A | \in \Real$, there is a
unique \nmulti\ having highest weight. Conversely, any \nmulti\
corresponds to a permutation on $A$.

\begin{remark}\label{rmk:graph} Given a
digraph $\grph$ where $\ideg(v) \geq 1$ (resp. $\odeg(v) \geq 1$),
for each $v \in V$, then $G$ contains a simple cycle. Indeed,
otherwise $G$ must have a source (resp. sink),  $v \in V$, in
contradiction to $\ideg(v) \geq 1$ (resp. $\odeg(v) \geq 1$),
respectively.
\end{remark}

%For any graph $\grph = (E,V)$ and any set $W \subset V$ of
%vertices, we define $E_W$ to be the set of edges adjacent to
%vertices in $W$; in other words $(v_1,v_2) \in E_W$ if $v_1 \in W$
%or $v_2 \in  W.$ Then we define the subgraph $\graph_W =
%(E\setminus E_W,V \setminus W)$.

In the following exposition we write $A \preceq 0^\nu$ for a
matrix $A$ all of whose  entries are $\preceq 0^\nu$ and assume
$\Det{A} \neq \minf$.

\begin{lemma}\label{lem:singleRegualrZero} Assume $A \preceq 0^\nu
$ is an $n\times n$ matrix, each of whose columns (resp. rows)
contains at least one  $0$-entry or  $0^\nu$-entry, then, for some
$i$, $a_\si{i}{\msig (i)} \in \{0, 0^\nu \}$.
\end{lemma}

\begin{proof}
We may  assume $\msig \in S_n$ is the identity. Suppose $\msig$
does not involve any $0$-entry or  $0^\nu$-entry, and let $G'_A$
be the reduced graph of $A$ obtained by erasing all edges of $G_A$
having weights $\prec 0$. Thus, $G'_A$ has no self loops and
$\oDeg(v_i) \geq 1$ for each $v_i \in V$. But, by Remark
\ref{rmk:graph}, $G'_{A} $ has a simple cyclic $C = (v_{i_1,i_2},
\dots, v_{i_{k-1},i_k},v_{i_k,i_1})$ of wight $0$ or $0^\nu$ which
contradicts the maximality of $|A|$, since $w(v_{i_1,i_2}, \dots,
v_{i_{k-1},i_k},v_{i_k,i_1}) \succ w(v_\si{i_1}{i_1}, \dots,
v_\si{i_k}{i_k})$.
\end{proof}

\begin{corollary}\label{thm:singularitySingelTag}
 An $n \times n$ matrix $A \preceq 0^\nu$,
 each of whose columns (resp. rows) contains at least one
$\uuu{0}$-entry, is singular.
\end{corollary}

\begin{lemma}\label{thm:singularitySingelTag}
 An $n \times n$ matrix $A \preceq 0^\nu$,
 each of whose columns (resp. rows) contains either, at least two
$0$-entries or a $0^\nu$-entry, is singular.
\end{lemma}
\begin{proof} We may assume $\msig$ is the identity. If $a_\si{i}{i} =
0^\nu$,  for some $i$, we are done. Otherwise, let $G'_A$ be the
reduced graph of $A$ obtained by erasing all edges of $G_A$ having
weights $\prec 0$, i.e. $G'_A$ has only edges weighted $0$ or
$0^\nu$. Considering an edge with wight $0^\nu$ as a duplicated
edge, $\iDeg(v_i) \geq 2$, for each vertex of $G'_A$. Thus, by
Lemma \ref{lem:singleRegualrZero}, $G'_A$ has a self-loop.

Erase all the self-loops of $G'_A$ (in particular, non of these
self-loops is a duplicated edge), and denote the new graph by
$G''_A$. Each vertex of $G''_A$ has $\iDeg(v_i) \geq 1$ and thus,
by Remak \ref{rmk:graph}, $G''_A$ has a simple cycle, say $C =
(v_{i_1}, v_{i_2}, \dots , v_{i_k}, v_{i_1})$. This means that
$G'_A$ must have a self-loop for each $v_{i_u}$, $u = 1, \dots,
k$, since otherwise we would get a contradiction  to the
maximality of $\gm$, that is the $\nu$-evaluation of $\msig$ in
$A$. Thus, the permutation obtained from $\msig$ by replacing
these self-loops with the simple cycle $C$ has the same
$\nu$-evaluation $\gm$ as $\msig$ has.
\end{proof}

\begin{remark}\label{lem:changeMinor} Any $n \times n$ singular matrix $A$ has an $n \times (n-1)$
submatrix which can be replaced by $\pi(A_{i,j})$ without changing
the singularity of $A$. Indeed, given a permutation $\msig \in
S_n$, let $\gm$ denote the product $a_{1,\msig(1)} \cdots
a_{n,\msig(n)}$,  then:
\begin{enumerate}

    \item If  $\gm \in \etUnit$ then  there is some $a_{i,
    \msig(i)} \in \etUnit$ and we can replace all the columns  $j \neq \msig(i)$.
    (Note that if $\gm = \minf$ then $\Det{A} = \minf$.)
    \item When $\gm \in \Real$, since $A$ is singular, there are at
    least two different permutations $\msig_1$ and $\msig_2$, so we can replace
    any possible $n \times (n-1)$
submatrix  of  $A$.

\end{enumerate}
\end{remark}

%******************************* SECT  *********************************
\section{The rank of tropical matrices}

The notion of linear dependence, cf. Definition
\ref{def:tropicalDependent}, provides a natural definition for the
rank of a tropical matrix:
\begin{definition}\label{def:matrixRank}
The \textbf{tropical rank}, denoted $\rnk(\tA)$, of a matrix $\tA$
is defined to be the maximal number of tropically independent rows
in $\tA$.
\end{definition}
 By this definition a nonzero matrix, i.e.  $A \neq (\minf)$, can
have rank $0$; for example any matrix all of whose entries are
ghost has rank $0$.

 The following familiar properties of matrix
rank are easily checked for our tropical rank:
\begin{enumerate} \eroman

    \item  The rank of a submatrix can not exceed the rank of the whole
    matrix. \pSkip

    \item The rank is invariant under reordering of either rows or
    columns. \pSkip

    \item The rank is invariant under (tropical) multiplication of rows or
columns by real constants, and under insertion of a row or a
column obtained as a combination of others. \pSkip
\end{enumerate}

\noindent The last property is true, since otherwise if it would
changed the rank, then one could a priori choose this combination
to obtain a lower  rank. On the other hand, just as in the
classical theory, vectors which are tropically dependent in the
initial collection are still tropically dependent in any extended
collection.

Later, we prove that $\rnk(A)$ is also equal to the maximal number
of its independent columns, and therefore that the tropical rank
of a matrix and its transpose are the same, cf. Corollary
\ref{thm:rankOfTranspose}.

%******************************* SECT  *********************************
%\secSpc
\subsection{Tropical regularity and tropical dependence}
The theorems in this subsection provide the connection between
tropical dependence and tropical regularity, we open with the
special case of matrices having determinant equals  $\minf$.

\begin{definition}\label{rdef} We say that a set $r_1, \dots, r_m$ of rows has
\textbf{rank defect} $k$ if there are $k$ columns, which we denote
as $c_1, \dots, c_k$, such that $a_{i,j} = \minf$ for all $1 \le i
\le m$ and $1\le j \le k$.
\end{definition}

For example, the rows  $(2,\minf,2,\minf), (\minf,\minf,\minf,2),
(1,\minf,\minf,\minf)$ have rank defect 1, since they are all
$\minf$ in the second column.

\begin{proposition}\label{ssing} Given an $n\times n$ matrix $A$, then  $\Det{A} = \minf$,
iff for some $1 \le k \le n$, $A$ has $k$ rows having rank defect
$n+1-k.$\end{proposition}
\begin{proof} $(\Leftarrow)$ If $k=n$ then this is obvious, since
some column is entirely $\minf$. If $n>k$, we take one of the
columns $c_j$ other than $c_1, \dots, c_k$ of
Definition~\ref{rdef}. Then for each $i$, the $(i,j)$ minor
$A_{i,j}$ has
  at least $k-1$ rows with rank defect
$(n-1)+1-k$, so has determinant $\minf$ by induction; hence $\Det{
A } = \minf,$ by Formula~\eqref{def:tropicalDetByMinors}.

$(\Rightarrow)$ We are done if all entries of $A$ are $\minf,$ so
assume for convenience that $a_{n,n}\ne \minf$. Then $
\Det{A_{n,n}} = \minf$, so, by induction, $A_{n,n}$ has $k\ge 1$
rows of rank defect $(n-1)+1-k = n-k$.
We may assume that $a_{i,j} = \minf$ for $1 \le i \le k$ and $1
\le j \le n-k.$ Thus, we can partition $A$
as the matrix $$ A = \( \begin{matrix} \minf & B' \\
B'' & C
\end{matrix}\) ,$$
where $\minf$ denotes the $k \times n-k$ matrix all of whose
entries $\minf$, $B'$ is a $k\times k$ matrix, $B''$ is an $n-k
\times n-k$ matrix, and $C$ is an $n-k \times k$ matrix.
Accordingly, $\Det{B' }= \minf$ or $\Det{B'' } = \minf$.

If $\Det{B' } = \minf$,
 then, by induction, $B'$ has $k'$ rows of rank defect
$k+1-k',$ thus, the same $k'$ rows in $A$ have rank defect $(n-k)
+ k+1-k' = n+1-k',$ and we are done taking $k'$ instead of $k$.
If $\Det{B'' } = \minf$,
 then, by induction, $B''$ has $k''$ rows of rank defect
$(n-k)+1 -k'',$  these $k+k''$ rows in $A$ have rank defect
$n+1-(k+k''),$ and we are done, taking $k+k''$ instead of $k$.
 \end{proof}

\begin{theorem}\label{thm:IndependentToRegular}
An $n \times n $ tropical matrix  of rank $< n$ is tropically
singular.
\end{theorem}

\begin{proof}
Let $r_i$ denotes the $i$'th row of $A$.  Since $\rnk(A) < n$,
there are $\al_1, \dots, \al_n \in \eReal$, not all of them
$\minf$, such that $\al_1 r_1 \TrS \cdots \TrS \al_n r_n \in
\etUnit^{(n)}$.

If $\al_i = \minf$ for some $i$, say $i = n$, then the first
$(n-1)$ rows are tropically dependent and each minor $A_\si{n}{j}$
has rank $< (n-1)$. But then,  by induction, each $A_\si{n}{j}$ is
singular and, by Formula \eqref{def:tropicalDetByMinors}, $\Det{ A
} = \bigoplus_j a_\si{n}{j}\Det{A_\si{n}{j}} \in \etUnit$.

Assuming all $\al_i$'s are in $\Real$, we replace each row $r_i$
of $A$  by $\al_i r_i$ and have $ \bigoplus_i r_i \in
\etUnit^{(n)}. $ Let $\bt_i$ denotes the maximal value in each
column $j$. If $\bt_i = b_i^\nu$ we take $b_i$ instead of $\bt_i$;
when $\bt_j = \minf$, for some $j$, we replace it by an arbitrary
real. Let $A'$ be the matrix obtained by dividing each column $j$
of $A$ by $\bt_j$, accordingly $A'$ satisfies the conditions of
Lemma \ref{thm:singularitySingelTag} and is singular. Since
regularity/singularity is preserved under the above operations,
$A$ is singular.
\end{proof}

\begin{example}
Consider a $2 \times 2 $ matrix $(a_\si{i}{j})$ with rank $= 1$.
Then, there are $\al_1, \al_2 \in \eReal$ such that
$\al_1(\aaa_\si{1}{1},\aaa_\si{1}{2}) \TrS \al_2
(\aaa_\si{2}{1},\aaa_\si{2}{2}) \in \etUnit^{(2)}$.  Note that
$\al_i \neq \tUniS$,  since otherwise for $k \neq i$
$(\aaa_\si{k}{1},\aaa_\si{k}{2}) \in \etUnit^{(2)}$ which would
contradict the data that $\rnk(A) = 1$. Replacing each row $r_i$
of $A$ by $\al_i r_i$ and expanding the determinant we get
$$
 (\al_ 1\aaa_\si{1}{1}) (\al_2 \aaa_\si{2}{2}) = (\al_2
\aaa_\si{1}{2})(\al_2 \aaa_\si{2}{1}) \Dir \al_1 \al_2 |\tA| =
\al_1 \al_2( \aaa_{1,1}\aaa_{2,2} \TrS  \aaa_{1,2}\aaa_{2,1}) \in
\etUnit,
$$
i.e. $\tA$ is tropically singular.

\end{example}

\begin{theorem}\label{thm:regularityToIndependent}
An $n \times n$  matrix $A$ has rank $< n$ iff $\tA$ is tropically
singular.
\end{theorem}

\begin{proof}
\noindent ${\bf (\dir )}$ By Theorem
\ref{thm:IndependentToRegular}.

${ \bf (\Leftarrow)}$ Assuming that $\tA$ is singular we need to
prove that the rows of $A$ are tropically dependent. Since parts
of the proof is  by induction $n$, the size of $A$, we assume the
theorem is true for $(n-1)$; the case of $n =1$ obvious. (The case
of $n=2$ is provided in Example \ref{exm:1}.)

Throughout this prove, we assume $\msig$ is the identity i.e.
$$ (\Det{A})^\nu = \gm^\nu =  (a_{1,1} \cdots a_{n,n})^\nu,$$
this hypothesis is not affected by multiplying through any row or
column by a given $\al \in \Real$. We also remark that when
determining the dependence coefficients $\al_i$'s, we may assume
the relevant $a_{i,j}$ are  in $\eReal$, since otherwise for
$a_{i,j} = b^\nu$ we take $b$ instead.

\textbf{Case I:} For notational convenience, if $A$ has an $m
\times m$ singular submatrix $A'$ with $\pi(|A'|) = \pi(
a_{i_1,i_1} \cdots a_{i_m,i_m})$, renumbering the indices, we
assume that the singular submatrix $A'$ with the minimal $m$ is
the upper left submatrix of $A$, in particular if $a_{i,i} \in
\etUnit$, for some $i$, renumbering the indices we may assume
$a_{1,1} \in \etUnit$.

Let
\begin{equation}\label{eq:al}
\al_i = \pi\(|A_{i,1}|\),
\end{equation}
excluding the case when all $\al_i$'s are $\minf$, see Case II, we
claim that $\bigoplus_i \al_i r_i \in \etUnit^{(n)}$, i.e.
\begin{equation}\label{eq:toProv} \bigoplus_i \al_i a_{i,j} \in
\etUnit, \quad \text{ for each $j = 1,\dots,n$}.
\end{equation}

Suppose $j = 1$, then $\bigoplus_i \al_i a_{i,1} \in \etUnit$,
since this is just the expansion of $| A|$ along the first column
of $ A$, i.e. $\bigoplus_i \al_i a_{i,1} = |A|$. Indeed, if $m
=1$, i.e. $a_{1,1} \in \etUnit$, we are done.  Otherwise,
$(a_{1,1} |A'_{1,1}|)^\nu = (a_{1,i} |A'_{1,i}|)^\nu = |A'|$ for
some $1 < i \leq m $.  Thus, since $\gm  = a_{1,1} |A'_{1,1}| \bt
= a_{1,i} |A'_{1,i}| \bt$, up to $\nu$, we have   $\al_1 a_{1,1} =
\al_i a_{1,i}$.

Assume $ \al_{\ell} a_{\ell,j}$, with $j > 1$,  is a component
with maximum $\nu$-value in the sum $  \al_{1} a_{1,j}\TrS \cdots
\TrS \al_{n} a_{n,j} $. If  $ \al_{\ell} a_{\ell,j} = \minf$ we
are done, otherwise $\al_\ell$ and $a_{\ell,j}$ are not $\minf$,
then
\begin{equation*}\label{eq:s1}
\al_{\ell}  = \pi \( |A_{\ell,1}|\) = \pi\(a_{1,\sig(1)} \cdots
a_{\ell-1,\sig(\ell-1)} a_{\ell+ 1,\sig(\ell+1)}  \cdots
a_{n,\sig(n)} \) , \qquad \sig(i) \neq 1,
\end{equation*}
for some $\sig \in S_n$. Let  $u $ be the index for which $\sig(u)
= j$, i.e. $u \neq \ell$, then
$$ \frac{|A_{\ell,1}|}{a_{u,j}} =   \frac{\(a_{1,\sig(1)} \cdots
a_{u,j}   \cdots a_{n,\sig(n)} \)}{ a_{u,j}} = \(a_{1,\sig(1)}
\cdots a_{u-1,\sig(u-1)} a_{u+1,\sig(u+1)}\cdots
a_{n,\sig(n)}a_{\ell,j} \),
$$
 up to $\nu$, which must be equal to
 $$\frac{\Det{A_{u,1}}}{a_{\ell,j}} = \frac{\(a_{1,\sig'(1)}
\cdots a_{u-1,\sig'(u-1)} a_{u+1,\sig'(u+1)}\cdots a_{\ell,
\sig'(\ell)} \cdots a_{n,\sig'(n)}a_{\ell,j} \)}{a_{\ell,j}}, $$
 since
otherwise $ a_{u,j} |A_{u,1} | \succ  a_{\ell,j}|A_{\ell,1} |$,
contrary to hypothesis. So, $a_{u,j} |A_{u,1} |$ and $
a_{\ell,j}|A_{\ell,1} |$ are two different terms in Formula
\eqref{eq:toProv} having  a same $\nu$-value.

\textbf{Case II:} When $\Det{A } = \minf $, with all $\al_i =
\Det{A_{i,1}}$ are $\minf$, we take $m$ maximal such that $A$ has
an $m\times m$ submatrix of determinant $\ne \minf $, and let
$\gm$ denote the determinant of the $m\times m$ submatrix $A_m$ of
$A$ of maximal $\nu$-value. By induction, we may assume that $m=
n-1.$ Furthermore, it is enough to find a dependence among the $k$
rows obtained in Proposition~\ref{ssing}, so, again, by induction,
we may assume that $k = n$, and the entries in the first column
are all $\minf $. Since $a_{1,1} = \minf$ and $\Det{A_m} \neq
\minf$, namely $A$ has an $(n-1)\times(n-1)$  minor whose
determinant $\neq \minf$, the proof is then completed by the same
arguments of Case I.
\end{proof}

\begin{corollary}\label{thm:rankN}
A matrix $\tA\in \matS{n}{n}(\Trop)$ has rank $n$ iff $\tA$ is
non-singular iff $\tA$ is pseudo inevitable.
\end{corollary}
\begin{proof} The proof is derived from Theorem
\ref{thm:regularityToIndependent} and Theorem
\ref{thm:inverseMatrix}.
\end{proof}
This corollary  provides the complete tropical analogues to the
well known classical relations between regularity, invertibility,
and rank of matrices.

\begin{example}\label{exm:1}
Suppose $\tA = (a_{i,j})$ is a $2 \times 2$ singular matrix, i.e.
$|\tA| = \aaa_{1,1}\aaa_{2,2} \TrS \aaa_{1,2}\aaa_{2,1}$. If $|A|
= \minf$ then $A$ has a $\minf$ row, say $r_1$, then set $\al_2 =
\minf$  and take arbitrary real $\al_1$. Otherwise, $\msig$ is the
identity, so take
$$ \al_1 =  \pi(| A_{1,1}|) = \pi(a_{2,2}), \qquad \text{and}
\qquad \al_2 =  \pi(| A_{2,1}|) = \pi(a_{1,2}).$$
Then, $$\al_1 r_1 \TrS \al_1 r_1 = \pi(a_{2,2})(a_{1,1}, a_{1,2})
\TrS \pi(a_{1,2})(a_{2,1}, a_{2,2}) \in \etUnit^{(2)}.$$
\end{example}

\begin{example}
Consider the  matrix
$$ A =
   \vvMat{1}{4}{-1}
        {1}{0}{6}
        {-4}{1}{3}
%\overset{\rw'_2 := (0)\rw_1 \TrS \rw_2}{\Dir}
%   \vvMat{1}{4}{-1}
%        {\uuu{1}}{4}{6}
%        {-4}{1}{3}
%\overset{\rw'_3 := (-3)\rw'_2 \TrS \rw_3}{\Dir}
%   \vvMat{1}{4}{-1}
%        {\uuu{1}}{4}{6}
%        {\uuu{-2}}{\uuu{1}}{\uuu{3}},
        $$
whose determinant equal $\uuu{8}$, and thus is singular. The
tropical dependence of the rows of $A$ is given by $\al_1 =
|A_{1,1}| = 7$, $\al_2 = |A_{2,1}| = 7$, and $\al_3 = |A_{3,1}| =
 10$.
\end{example}

\begin{theorem}\label{thm:dependenyOfVectors}
Any $k > n$ tropical vectors in $\Trop^{(n)}$  are tropically
dependent.
\end{theorem}
\begin{proof} Assume $v_1, \dots, v_{n+1}$ are independent
vectors in $\Trop^{(n)}$ and consider the $(n+1) \times n$ whose
rows are these vectors. Extend this matrix by duplicating one of
the columns to get a singular matrix, cf. \cite[Theorem
2.5]{zur05TropicalAlgebra}, whose rows are tropically dependent by
Theorem \ref{thm:regularityToIndependent}, a contradiction.
\end{proof}

Our next goal is to show that the rank of an $m \times n$ matrix
is determined as the maximal size of whose maximal nonsingular
minor, rather than  by a collections of minors of smaller sizes.

\begin{theorem}\label{thm:maximalMinor}
An $m \times n$ matrix $A$ of rank $m$, $m \le n$, has an $m
\times m$ nonsingular minor $A_{\max}$.
\end{theorem}

\begin{proof} The case when $m=n$ is obvious by Theorem
\ref{thm:regularityToIndependent}, we proceed by induction on $n$.
Let $A'$ denote the $m \times (n-1)$ submatrix of $A$ obtained by
erasing the last column and let $A''$ be the submatrix of $A$
obtained by erasing the first column. Assuming both $A'$ and $A''$
have rank $ < m$, we aim for a contradiction.

Throughout this prove, to make the exposition clearer, we often
use matrix products to describe sums; for example we write $(a_1,
\dots, a_n) (b_1, \dots, b_n)^\trn$ for the sum $\bigoplus_i a_i
b_i$.

Denoting the rows of $A'$ as $r'_i$, since $\rnk(A') < m$, there
are $\al'_1, \dots, \al'_m \in \eReal$, not all of them $\minf$,
such that
$$ \al'_1 r'_1 \TrS \cdots \TrS \al'_m r'_m \in \etUnit^{(n-1)}.$$
We write $\bal'$ for the $m$-tuple $(\al'_1, \dots, \al'_m)$ and
define  $\bal''$ by the same way for $A''$.

We show that there are $\mu', \mu'' \in \eReal$, for which $\bbt =
\mu' \bal' \TrS \mu'' \bal''$ determines a dependence on $A$. We
also need to verify that each entry of $\bbt$ is in $\eReal$.

Let $r_i$ denote the $i$'th row of $A$, $c_j$  denote the $j$'th
column of $A$,  and write
\begin{equation}\label{eq:mprod}
    (\mu', \mu'') \(\begin{array}{ccc}
  - &  \bal' & - \\
  - &  \bal'' & - \\
\end{array}%
\)\(\begin{array}{c|ccc|c}
    a_{1,1} & \cdots & a_{1,j} & \cdots & a_{1,n} \\
    \vdots & & \vdots & & \vdots \\
    a_{m,1} & \cdots & a_{m,j} & \cdots & a_{m,n} \\
  \end{array}\) = (b_1,\dots,b_n).
\end{equation}
Since  $c_j$, for  $ j = 2, \dots, n-1$, is a column of both $A'$
and $A''$, and thus $(\bal')(c_j) \in \etUnit$ and $(\bal'')(c_j)
\in \etUnit$, it is clear that $b_j$ is ghost for each $ j = 2,
\dots, n-1$. So,  by leaving only the first and the last column of
$A$, we reduce Formula \eqref{eq:mprod} and write
\begin{equation}\label{eq:mprod2}
    (\mu', \mu'') \(\begin{array}{ccc}
  - &  \bal' & - \\
  - &  \bal'' & - \\
\end{array}%
\)\(\begin{array}{cc}
    a_{1,1}  & a_{1,n} \\
    \vdots   & \vdots \\
    a_{m,1} & a_{m,n} \\
  \end{array}\) = (\mu', \mu')
\begin{array}{c}  = B \\
  \(%
\begin{array}{cc} \hline
  (\bal')(c_1) & (\bal')(c_n) \\
  (\bal'')(c_1) & (\bal'')(c_n) \\
\end{array}%
\).\\
\end{array}
\end{equation}

It easy to see that $B$ is singular, just expand $\Det{B}$ to get
$$ \Det{B} = (\bal')(c_1)(\bal'')(c_2) \TrS  (\bal'')(c_1)(\bal')(c_n) =
%(\bigoplus_i \al'_i a_{i,1}) (\bigoplus_j \al''_j a_{j,n}) \TrS
%(\bigoplus_i \al''_i a_{i,1}) (\bigoplus_j \al'_j a_{j,n})=
\bigoplus_{i,j} \al'_i a_{i,1} \al''_j a_{j,n} \TrS
\bigoplus_{i,j} \al''_i a_{i,1} \al'_j a_{j,n},$$ and thus
$\Det{B} \in \etUnit$. Therefore, by Theorem
\ref{thm:regularityToIndependent}, $B$ has rank $< 2$ and whose
rows are tropically dependent;  this means that there are $\mu'$,
$\mu''$ in $\Real$ for which $(\mu',\mu'') B \in \etUnit^{(2)}$.

Assume there is $i$ for which $\bt_k a_{1,i} \succ
\bigoplus_{h\neq i } \bt_h a_{1,h}$, where $\bt_i \in \etUnit$.
But, $\bt_i = \mu'\al'_i + \mu''\al''_i $ and, by hypothesis on
$A'$, there is $\al'_k$, $k \neq i$, for which  $\al'_k a_{1,k} =
\al'_i a_{1,i}$, equivalently $\mu' \al'_k a_{1,k} = \mu'\al'_i
a_{1,i}$. Therefore, $\bt_i$ can be replaced by $\pi(\bt_i)$. By
interchanging $\al'_k$ by $\al''_k$, and taking all indices with
respect to $A''$, the same argument is applied to the column
$c_n$.

This shows that $\pi(\bbt)$ determines a tropical dependence on
the rows of $A$, a contradiction to the data $\rnk(A) = m$. Thus,
either $A'$ or $A''$ has rank $m$, and by the induction hypothesis
has an $m \times m$ nonsingular minor, also a minor of $A$.
\end{proof}

\begin{corollary}\label{thm:rankToMinors}
An $m \times n$ matrix $\tA$ has rank $k$ iff its maximal
nonsingular minor is of size $k \times k$.
\end{corollary}
\begin{proof}  $\tA$ can not have a minor $\tA_K$  of rank grater
than $k$, since otherwise $\rnk(\tA)$ would be grater than $k$.
The proof is then completed by Theorem \ref{thm:maximalMinor}
applied to $\tA_K$.
\end{proof}

\begin{corollary}\label{thm:rankOfTranspose}
The rank of a matrix and the rank of its transpose are the same.
\end{corollary}

\begin{proof}
The rank of $A$ and $A^\trn$ are both equal to the size of the
maximal nonsingular minor.
\end{proof}

\begin{corollary}\label{thm:columnsDependency}
The rank of a matrix is equal to size of a maximal independent
subset of its columns.
\end{corollary}
%******************************* SECT  *********************************
\secSpc \subsection{Relations to former settings }

Recall that for real matrices our definition of singularity
coincides with the known definition for matrices over
$\eMaxPlusAlg$, cf Remark \ref{rmk:regularity}. In
\cite{Develin2003}, Develin, Santos, and Sturmfels, define the
tropical rank of an $n\times n$ matrix $\tA$ over $\eMaxPlusAlg$
to be the largest integer $k$ such that $\tA$ has a $k \times k$
nonsingular minor, we denote this type of rank by $\rnk_D(\tA)$
and the corresponding nonsingular minor of maximal size by
$A_{\max}$. To emphasize, $\rnk_D(\tA)$ is given only for matrices
with real entries without any notion of linear dependence. Our
work bring in the notion of linear dependence, and in the light of
Corollary \ref{thm:rankToMinors} we have:
\begin{proposition}\label{thm:DevelinRank} When $\tA$ is a  real matrix, i.e.
$A \in \matS{n}{n}(\eReal)$, the tropical rank as in Definition
\ref{def:matrixRank} coincides with that of Develin, Santos and
Sturmfels, i.e.
$$\rnk_D(\tA) = \Size(A_{\max}) = \rnk(\tA),$$
where $A_{\max}$ is a nonsingular minor of maximal size of $A$.
\end{proposition}
\begin{proof} Immediate by  Corollary
\ref{thm:rankN} and Corollary \ref{thm:rankToMinors}.
\end{proof}

Therefore, concerning real matrices, our rank preserves also the
known relation to Barvinok and Kapranov ranks \cite{Develin2003},
denoted respectively as $\rnk_B(\tA)$  and $\rnk_K(\tA)$,  that is
$$ \rnk(\tA) \ \leq \ \rnk_K(\tA) \ \leq \ \rnk_B(\tA),$$
for each  $A \in \matS{n}{n}(\eReal)$.

\begin{remark}
The computation of each of the above ranks for matrices over
$\eMaxPlusAlg$ has been proven to be $NP$-complete
\cite{Kim:130853}. Thus, in the view of  Proposition
\ref{thm:DevelinRank}, computating our rank is $NP$-complete as
well.
\end{remark}

The below definition and proposition were introduced and proven in
\cite{Butkovic1985}, are applied only to square matrices defined
over ``pure reals", i.e. non of the matrix entries is $\minf$.
\begin{definition}\label{def:stronglyIndepentent}
The columns of an ${n \times n}$ (pure real) matrix $\tA$ are
\textbf{strongly linearly independent} if there is a column vector
$\tV \in \Real^{(n)} $ such that the tropical linear system $\tA
\tX = \tV$ has a unique solution $\tX \in \Real^{(n)}$. A square
matrix is \textbf{strongly regular} if its columns are strongly
linearly independent.
\end{definition}
\begin{proposition}\label{thm:StronglyRegular} For a
square (pure real) matrix,  strongly regular and tropically
nonsingular are equivalent.
\end{proposition}
Note that according to Definition \ref{def:stronglyIndepentent},
strongly independence is based on both the existence of a solution
and on its uniqueness. Using this notion of dependence, the
development becomes very difficult and not intuitive.
\begin{corollary}\label{thm:StronglyRegularDependency}
The tropical rank of a pure real matrix equals the largest size of
a strongly linearly independent subset of its columns.
\end{corollary}
\begin{proof}
Immediate, by Proposition \ref{thm:StronglyRegular} and the
equality of $\rnk_D(\tA) = \rnk(\tA)$ for real matrices.
\end{proof}

%******************************* SECT  *********************************
\secSpc
\section{Solutions of homogeneous linear systems}

Recall that the zero set of a polynomial $f \in \Trop[\lm_1, \dots
\lm_n]$ in $n$ indeterminates $\lm_1, \dots, \lm_n$ is defined as
$$ Z(f) =\{ \bfa \in \Trop^{(n)} \ | \  f(\bfa) \in \etUnit\},$$
where $\bfa$ stands for $(a_1,\dots, a_n)$, cf. \cite[Defintion
2.2]{zur05SetNullstellensatz}. A polynomial $f$ is said to be
homogenous if all of whose monomials are of the same degree.

Using this notion of zeros, we say that a system $S$ of $m$
homogenous linear equations
\begin{equation}\label{eq:linearSystem}
   \begin{array}{ccccccc}
   f_1 & = & \aaa_{1,1} \lm_1 &\TrS &\cdots &  \TrS  & \aaa_{1,n} \lm_n,\\
      \vdots &   & \vdots & & & & \vdots \\
    f_m & = &  \aaa_{m,1} \lm_1 & \TrS & \cdots &  \TrS & \aaa_{m,n} \lm_n,  \\
\end{array}
\end{equation}
has a solution if all equations have a common zero, i.e. there
exist $\bfa \in \Trop^{(n)}$ such that $f_i(\bfa) \in \etUnit$,
for all $i = 1, \dots,m$. When $\bfa \in \Real^{(n)}$ we say that
the solution is pure real.

\begin{remark} Any $\bfa \in \etUnit^{(n)}$ is also a solution of
a system \eqref{eq:linearSystem}. We can also also have mixed
solutions, these are solutions for which $\bfa$ has entries in
both $\Real$ and $\etUnit$.
\end{remark}

As usual a system $S$ of the form \eqref{eq:linearSystem} can be
written in matrix terms as $ A_S \Lm^\trn $, where $A_S$ is the $m
\times n$ coefficients matrix of the system $S$ and $\Lm = (\lm_1,
\dots, \lm_n)$.

\begin{theorem}
A system $S$ of $n$ homogenous linear equations has a pure real
solution iff  the corresponding coefficient matrix $\tA_S$ is
singular.
\end{theorem}
\begin{proof}  Obvious  by Corollary
\ref{thm:columnsDependency}. \end{proof}

Note the a system $S$ with $A_S$ nonsingular can also have mixed
solutions.

%******************************* Reference *********************************
\small
%\bibliographystyle{abbrv}
%\bibliography{../../../../bib/dfz}

\end{document}